\input amstex
\documentstyle{amsppt}
\magnification=1200
\input prepictex
\input pictex
\input postpictex
\input diagrams.sty

\pagewidth{6.5 true in}
\pageheight{9.0 true in}
\NoBlackBoxes
\topmatter
\title
Heights of Ideals of Minors
\endtitle
\rightheadtext{Heights of Ideals of Minors}
\author
David Eisenbud, Craig Huneke and Bernd Ulrich
\endauthor
\address
Mathematical Sciences Research Institute,
1000 Centennial Dr.,
Berkeley, CA 94720
\endaddress
\email
de\@msri.org
\endemail
\address
Department of Mathematics, University of Kansas,
Lawrence, KS 66045
\endaddress
\email
huneke\@math.ukans.edu
\endemail
\address
Department of Mathematics,
Purdue University,
West Lafayette, IN 47907
\endaddress
\email
ulrich\@math.purdue.edu
\endemail
\thanks
All three authors were partially supported by the NSF.
\endthanks

%%%%%%%%%%%%%%%%%%%%%%%%%%%%%%%%%%%%%%%%%%%%%%%%%%%%%%%%%%%%%%%%%%%

\def \fix#1 {{\hfill\break \bf (( #1 ))\hfill\break}}

\def\A{\frak A}
\def\E{\Cal E}
\def\F{\Cal F}
\def\R{\Cal R}
\def\B{\frak b}

\def\P{\Bbb P}
\def\Hom{\text{Hom}}

\def\bight{\text{bight}}

\def\n{\bold n}

\def\ra{\rightarrow}
\def\lra{\longrightarrow}

\def\rank{\text{rank}}

\def\l{\ell}
\def\im{{\frak m}}

\def\fin{{\frak n}}

\def\p{\phi}
\def\var{\varphi}

\def\L{\Lambda}

\def\d{\delta}
\def\e{\epsilon}

\def\a{\alpha}
\def\b{\beta}

\def\B{{\frak b}}
\def\A{{\frak a}}

\def\Hom{\text{Hom}}

\def\hgt{\text{ht}}

\def\height{\text{height\;}}
%\def\height{{\rm height\;}}

%\rightline{version of 5/6/2001}\medskip
%%%%%%%%%%%%%%%%%%%%%%%%%%%%%%%%%%%%%%%

\date{May 1, 2002}
\enddate

\endtopmatter

\document

{\narrower
\noindent{\bf Abstract.} We prove new height inequalities
for determinantal ideals in a regular local ring, or more
generally in a local ring of given embedding codimension.
Our theorems extend and sharpen results
of Faltings [F] 
and Bruns [B1]. \bigskip
}

\centerline{\bf Introduction\rm}
\bigskip

Let $\varphi$ be a map of vector bundles on a variety $X$.  A
well-known theorem of Eagon and Northcott [EN] gives an upper bound
for the codimension of the locus where $\varphi$ has rank $\leq s$ for
any integer $s$.  

Bruns [B1] improved this result by taking into
account the generic rank $r$ of $\varphi$.  We shall see below that
unlike the Eagon-Northcott estimate, Bruns' Theorem is sharp only when
$X$ is singular.  The first goal of this paper is to give stronger
results when $X$ is nonsingular, and a little more generally.

Strengthening the Eagon-Northcott estimate in a different way from Bruns,
Faltings [F1] gave an improved bound for the case $s=r-1$ under the
additional assumption that $X$ is nonsingular and the cokernel of
$\varphi$ is torsion free. We also improve Faltings theorem to a
result valid for all $s$.

Our results are actually local. Let $R$ be a local ring,
and let $\varphi: R^m\to R^n$ be a matrix of rank $r$.
We write $I_i=I_i(\varphi)$ for the ideal generated by $i\times i$ minors
of $\varphi$, and we assume that $I_i\neq R$.
Bruns' Theorem says that
$$
\height\!(I_i) \leq (r-i+1)(m+n-r-i+1).
$$
This formula is sharp for every $m,n,r,i$: take $\varphi$ to be the
image of the generic $n\times m$ matrix
$$
\Phi=(x_{ij}) \qquad 1\leq i\leq n,  \quad 1\leq j\leq m
$$
over the ring $R=k[\{x_{ij}\}]/I_{r+1}(\Phi)$. Note that
this ring is singular for $r>0$.

Henceforth
in this introduction we shall assume that
$R$ is a regular local ring.
Under this hypothesis we can improve Bruns' bound as follows:

\proclaim{Theorem A}  $\text{\rm height}(I_i) \leq (r-i+1)(\max(m,n)-i+1)+i-1.$
\endproclaim

\noindent Theorem A is a weak form of Corollary 3.6.1 below.

One should compare this result with the ``trivial'' case where the matrix 
$\var$ contains only $r$ nonzero rows (if $m\geq n$) or $r$ nonzero
columns (if $n\geq m$).  In this case the codimension of the 
ideal of $i \times i$ minors is given by the ``Eagon-Northcott'' formula
$$ \height\!(I_i) \le (r-i+1)(\max(m,n)-i+1),$$
which is an equality if the nonzero rows (columns) of $\var$ are generic.  This 
formula coincides with ours when $i = 1$.  Theorem A is also sharp for 
generic alternating $3\times 3$ matrices when $i = 2$.

A particularly interesting situation is that where the cokernel of
$\varphi$ is torsion free (or even a vector bundle on the punctured spectrum).
In this torsion free case Faltings improved Bruns' bound
(for the $r\times r$ minors only) and showed
$$
\height\!(I_r) \leq n.
$$
Generalizing this to arbitrary size minors, and allowing the
cokernel to be any module which is not the direct sum of a free
module and a nonzero torsion module (that is, excluding the
``trivial'' case described above) we show that:

\proclaim{ Theorem B}   $\text{\rm height}(I_i) \leq n + (r-i)(\max(m,n)-i+1).$
\endproclaim

\noindent Theorem B is a weak form of Corollary 3.6.2 below.

Theorem B is sharp in the case where $n=3, m \ge 3, r = i = 2$ and $\var$ 
is the generic alternating $3 \times 3$ matrix followed by a $3 \times 
(m-3)$ matrix of zeros.  If on the other hand $\var$ has one generic 
column, $r-1$ generic rows, and the rest of its entries 0, then
$$ \height\!(I_i) = (r-i)(m-i+1) + \min\{m- i + 1, n-r+1\}.$$
This actual value is close to the bound given  by Theorems A and B.
Some less degenerate examples are given in section 4.

We can also ask for a bound on the
height of one ideal of minors modulo the ideal of
minors of the next larger size. By symmetry we may assume that $m\geq n$. 
We prove

\proclaim{Theorem C}
$
\text{\rm height}(I_i/I_{i+1}) \leq  \max(m-i+1, n)+r-i.
$
\endproclaim

\noindent Theorem C is a weak form of Corollary 3.9 below.

This result is comparable to Theorems A and B (or their sharpenings)
in the case $i=r$; but it does not follow from these results in general
because $R/I_{i+1}$ is not regular. However, if we have
good information about the higher order minors of $\varphi$,  as
in the case where the cokernel of $\varphi$ is an ideal, then Theorem C gives
results on the height of $I_i(\varphi)$ that are better than those
coming from Theorems A and B. In this way we reprove a theorem
of Huneke [Hu] and extend it as follows:

\medskip

\proclaim{Corollary D} Let $I$ be an ideal of $R$ of height $g$ 
that is minimally generated by $n$ elements.

\text{\rm (a) (Huneke)} If $I$ is not a complete intersection, that is $n>g$,
then the locus of primes $P$ such that $I_P$ is not a complete intersection
has codimension $\leq n+2g-1$.

\text{\rm (b)} If $R/I$ is a Cohen-Macaulay domain and $n>g+1$
then the locus of primes $P$ such that $I_P$ cannot be generated by $g+1$
elements has codimension $\leq n+2g-1$.
\endproclaim

\noindent Corollary D is a weak form of Example 3.11 below.

Huneke's result (which is sharp, for example, in case $I$ is the
ideal of $2\times 2$ minors of a $2\times 3$ matrix) improves a
formula of Faltings [F] by 1. One should compare this to a famous
conjecture of Hartshorne [Ha1]  saying that if $I$ is the homogeneous
ideal of a smooth projective variety $X$  which is not a complete
intersection, then the singular locus of $X$ has codimension $\leq 3g+1$
in the ambient projective space.

Both Theorems A and C are direct consequences of our other
main theorem, which gives the bound on the codimension of the ideals
of minors of a matrix
$\bar{\var}$ over a ring $\bar{R} =R/J$
obtained by reducing $\var$ modulo $J$. We may assume that $m\geq n$.
We write $\bar{r}$ for the rank of $\bar\var$ and set
$\delta=r-\bar{r}$.

\medskip

\proclaim{Theorem E}  $\text{\rm height}(I_i(\bar\var))  \leq
(\bar{r} -i)(m-i+\delta+1)+\max(m-t+1,n)+\delta.$
\endproclaim

\noindent Theorem E is a weak form of Theorem 3.1.1 below.  As with 
Faltings' work, we do not need $R$ to be regular, but can give bounds in
terms of certain embedding codimensions.

\bigskip

We now describe the key ideas of our proofs.  To establish height bounds 
for ideals of minors it is helpful to identify as ``many'' row ideals of 
$\var$ as possible that have ``small'' height.  As it turns out, the 
behavior of $\var$ in this respect is determined by the analytic spread 
$\l$ of $M = \text{Coker} (\var)$ (see Section 1 for the definition of 
analytic spread).  If $\l$ has the maximal possible value $n$ then all row 
ideals of $\var$ have height at most $r$, and (under weak conditions) the 
converse holds as well.  Thus, whenever $\l < n$ there have to exist row 
ideals whose height exceeds $r$.  On the other hand we  prove in this case 
that after a flat local base change, at least $\l$ row ideals have height 
$\le r-n+\l < r$.  To paraphrase, if the analytic spread of $M$ is not 
maximal, then the behavior of the row ideals is more unbalanced, but not 
necessarily worse for our purposes.  This is the content of Theorem 2.2, 
the main technical result of the paper.  A complicated induction then 
completes the proofs of our formulas in Section 3. 

\bigskip

We finish this introduction with a list of open problems specifically
suggested by the results of this paper. Of course the biggest open
problem is the conjecture of Hartshorne mentioned above.

\bigskip

{\bf Problem 1:}  Let $\var$ be a symmetric $n \times n$ matrix of rank 
$r$, and suppose that 2 is invertible in $R$ (but not necessarily that $R$ 
is regular).  We conjecture that for $i\leq r$,
$$
\height\!(I_i) \le \pmatrix n-i+2 \\ 2 \endpmatrix - \pmatrix n-r+1 \\ 2
\endpmatrix. 
$$
In Section 5 we prove this conjecture for the cases $i=1$ and $i=n-1$ if $R$ is
regular.

If the conjecture is true, it is sharp, for example for the generic matrix, taken 
modulo the ideal of $(r+1) \times (r+1)$ minors.  This formula is the 
analogue of Bruns' bound for general matrices; it is computed as the 
difference between the heights of the ideals of $i \times i$ and the $(r+1) 
\times (r+1)$ minors of a generic symmetric matrix.  Notice that the 
conjecture fails in characteristic 2, as can be seen by taking $\var$ to be 
a generic alternating $3 \times 3$ matrix and $i=2$.

\bigskip

{\bf Problem 2:}  Are there better bounds than the ones of Theorems A and B 
if we assume that $\var$  has no (generalized) rows or columns of zeros?

\bigskip

{\bf Problem 3:}  Are there better bounds if $\var$ is a matrix of linear 
forms?

\bigskip

{\bf Problem 4:}  Find sharp bounds assuming the ranks are small.  For 
example, what about $I_2$ for a $4 \times 4$ matrix of rank 2?  Is the 
height bounded by 3?

\bigskip
\bigskip

%%%%%%%%%%%%%%%%%%%%%%%%%%%%%%%%%%%%%%%%%%%%%%%%%%%%%%

\bigskip
\bigskip
\centerline{\bf 1. Basic Results\rm}
\bigskip

In this section we fix our notation and review some basic facts, mainly 
about Rees algebras of modules, that will be used throughout.

Let $R$ be a Noetherian ring and $I$ an ideal of $R$.  We write $\hgt(I)$ 
for the height of $I$ and $\bight(I)$ for its {\it big height} , which is 
the maximum of the heights of minimal primes of $I$.  Let  $M$ a 
finitely generated $R$-module and $\var$ an $n$ by $m$ matrix with entries 
in $R$.  By the {\it $i^{th}$ row ideal} of $\var$ we mean the ideal 
generated by the entries of the $i^{th}$ row of $\var$, and the {\it rank} 
of $\var$ is the integer $r = \max\{i | I_i(\var) \neq 0\}$.  We say that 
$M$ has a {\it rank} and write rank$(M) = e$ if $M \otimes_R K$ is a free 
$K$-module of rank $e$, with $K$ denoting the total ring of quotients of 
$R$.  Notice that if $\var $ presents $M$ and $M$ has a rank, then $r + e = n$.

Let $\var$ be a matrix presenting $M$ and $\underline{T} = T_1, \dots, T_n$ 
a row of variables.  The row ideals of $\var$ are related to the symmetric 
algebra $\text{Sym}(M)$ of $M$ via the homogeneous presentation 
$\text{Sym}(M) \cong R[T_1, \dots, T_n] /I_1(\underline{T} \cdot \var)$ 
(see also [EHU1], where this fact has been exploited 
systematically).  Since the symmetric algebra fails to be equidimensional 
in general, we are lead to consider the {\it Rees algebra} $\R(M)$ of $M$ 
instead.  The general notion of Rees algebra has been introduced in 
[EHU2, 0.1].  In the present paper however we will restrict ourselves 
to considering modules that have a rank.  In this case $\R(M)$ is equal to 
$\text{Sym}(M)$ modulo $R$-torsion.  We say that $M$ is of {\it linear 
type} if the natural map from $\text{Sym}(M)$ to $\R(M)$ is an 
isomorphism.  If $R$  has dimension $d$ and $E$ has a rank $e$, then $\dim 
\R(M) = d + e$ (see, e.g., [SUV, 2.2]).  Suppose in addition that $R$ is 
equidimensional, universally catenary and local.  Under this assumption 
$\R(M)$ is equidimensional. Thus we may write $\R(M) \cong R[T_1, 
\dots, T_n]/ \B$ with  $\bight(\B) = \hgt(\B)$.  In fact $\hgt(\B) = r$, 
the rank of any matrix with $n$ rows that presents $M$.

If $U$ is a submodule of $M$, we say that $U$ is a {\it reduction} of $M$ 
or, equivalently, $M$ is {\it integral} over $U$ if the ring $\R(M)$ is 
integral over its subalgebra $R[U]$.  In case $R$ is local with residue 
field $k$, the {\it analytic spread} $\l(M)$ of $M$ is defined to be the 
Krull dimension $\dim \R(M) \otimes_R k$.  The two notions are related by 
the fact that $\l(M) = \min\{\mu(U)|U \text{ a reduction of 
}  M\}$  whenever $k$ is infinite (here $\mu(-)$ denotes minimal number of 
generators).  One always has $\rank(M) \le \l(M) \le \mu(M)$  (see, e.g. 
[SUV, 2.3]), and the last inequality is an equality if and only if $M$  has 
no proper reduction, at least in the case of an infinite residue field.

Before describing more refined estimates, we need to review the property 
$G_s$, where $s$ is an integer:  A module $M$ of rank $e$ is 
said  to  satisfy $G_s$ if $\mu(M_P) \le \dim R_P + e - 1$ for every prime 
ideal $P$ with $1 \le \dim R_P \le s-1$.  What makes the concepts of 
integral dependence and analytic spread play a central role in this paper 
is their relation to the height of certain colon ideals:

\bigskip

\proclaim{Theorem 1.1} \text{\rm ([R, 2.5], [EHU3, 1.2])}  Let $R$ 
be an equidimensional universally catenary Noetherian local ring, let $M$ 
be a finitely generated $R$-module having a rank $e$, and let $U$ be a 
submodule of $M$ with $\mu(U) \ge e$.  If
$$\text{\rm ht}(U :_R M) > \mu(U) - e+1$$
then $U$ is a reduction of 
$M$.
\endproclaim

\bigskip

This theorem yields the upper bound  $\l(M) \le \mu(U)$ when the hypothesis 
is satisfied.   Conversely, one has:

\proclaim{Proposition 1.2} \text{\rm ([EHU3, 3.7bis])}  Let $R$ be a 
Noetherian local ring with infinite residue field, let $M$ be a finitely 
generated torsion free $R$-module having a rank $e$, and assume that $M$ 
satisfies $G_{s+1}$.  If $\text{\rm ht}(U :_R M) \le \mu(U) - e+1$ for 
every submodule $U$ generated by $e + s - 1$ general linear combinations of 
generators of $M$, then $\l(M) \ge e+s.$
\endproclaim

\bigskip

In a more general setting one still has the following weaker bounds:

\proclaim{Proposition 1.3} \text{\rm ([SUV, 4.1])}  Let $R$ be a Noetherian 
local ring and let $M$ be a finitely generated $R$-module having a rank $e$.
If $M$ is not a direct sum of a torsion module and a free 
module, and $M_P$ is free for every prime ideal $P$ with depth $R_P \le 1$, 
then $\l(M) \ge e + 1$.
\endproclaim

\bigskip

If $R \ra S$ is a homomorphism of rings, $J^c$ will denote the contraction 
to $R$ of an $S$-ideal $J$, and $-_S$ will stand for the functor $- 
\otimes_R S$.  We will denote $\text{\rm Hom}_R(-,R)$ by $-^*$.  The {\it 
embedding codimension} ecodim$(R)$ of a Noetherian local ring $(R,\im)$ is 
defined as the difference $\mu(\im) - \dim R$; equivalently, writing 
$\hat{R} \cong S/J$ with $(S,\fin)$ a regular local ring and $J$ an 
$S$-ideal contained in $\fin^2$, one has ecodim$(R) = \hgt(J)$.

%%%%%%%%%%%%%%%%%%%%%%%%%%%%%%%%%%%%%%%%%%%%%%%%%%%%%%%%%%%%%%%%%

\bigskip
\bigskip
\centerline{\bf 2. Choosing Row Ideals of Small Height\rm}
\bigskip

Let $R$ be an equidimensional universally catenary
Noetherian local ring, and $M$ a finitely
generated $R$-module having a rank $e$ with $n = \mu(M), \ell = 
\ell(M)$.  Theorem 1.1 shows that if $\ell = n$, then every row ideal of 
any matrix
minimally presenting $M$ has height at most $r = n-e$.  According to 
Proposition 1.2, the converse holds in case $M$ satisfies $G_{r+1}$. Thus,
whenever $\ell < n$ there tend to exist row ideals of height strictly
greater than $r$. On the other hand, we will prove below that
it is possible in this case to find ``many'' row ideals whose height is
strictly less than $r$. More precisely, over a flat local extension
ring $S$ of $R$ there exists a matrix $\p$ minimally presenting $M_S$
such that at least $\ell$ row ideals of $\p$ have height at most
$r-n+\ell = \ell-e$.  These row ideals are constructed as defining ideals 
of Rees algebras of certain modules.   The local homomorphism $R\ra S$ has 
a complete
intersection closed fiber, but regularity may fail to pass from $R$
to $S$. This will require some extra care since the
height of ideals in $S$ may no longer be subadditive.

We begin by recording a weaker version of the above estimate, which has
the advantage that $S$ can be chosen to be a localization of a polynomial ring
over $R$. This theorem was inspired by a result of Evans and Griffith saying
that if $R$ is a universally catenary domain with algebraically
closed residue field and $N$ is a finitely generated nonfree
$R$-module of rank $r$ then there exists a minimal generator $x\in N$
with $\hgt(N^*(x))\leq r$ (\cite{EG2, 2.12}).

\medskip

\proclaim{Theorem 2.1} Let $R$ be an equidimensional universally catenary
Noetherian local ring and let $M$ be a finitely
generated $R$-module with rank $e$. Write $\ell = \ell(M)$ and
  $r = \mu(M) - e$. Then there exists a local homomorphism
$R\ra S$, with $S$ a localization of a polynomial
ring over $R$,  and a minimal presentation matrix of $M_S$ over $S$
that has $\ell$ row ideals of height at most $r$.
\endproclaim

This result is a special case of the next theorem.
Before stating the theorem we remark on some notation and terminology.
Let $M$ be a finitely generated module over a local Noetherian
ring $(R,\im)$, and assume that $n = \mu(M)$.
When we speak of a \it generic \rm generating set of $M$ in a local ring
$R'$ obtained by $R$ by a purely transcendental residue field extension
we mean the following: Let $X = (x_{ij})$ be a generic $n$ by $n$ matrix
over $R$. Fix a generating set $m_1,...,m_n$ of $M$, and 
let $v_i = \sum_{i=1}^n x_{ij}m_j$, and set $R' = R[x_{ij}]_{\im R[x_{ij}]}$.
Then $M_{R'}:=M\otimes_R R'$ is generated by $v_1,...,v_n$. Furthermore, if $\ell = \ell(M)$,
then any $\ell$ of the $v_j$ form a minimal reduction of $M_{R'}$.
(This can be shown by first proving that the correct number of
generic elements always give Noether normalizations for finitely generated
algebras over fields, and which is explicitly shown in [FUV, Thm. 7.3].)
\medskip

\proclaim{Theorem 2.2} Let $(R,\im)$ be an equidimensional universally
catenary Noetherian local ring, let $M$ be a finitely generated $R$-module
with rank $e$, and set $n = \mu(M), \ell = \ell(M), r = n-e$. Let
$v_1,...,v_n\in M_{R'}$ be a generic generating set defined over
a local ring $R'$ that is obtained from $R$ by a purely transcendental
residue field extension and let $\psi$ be an $n \times m$ matrix presenting
$M_{R'}$ with respect to $v_1,...,v_n$. Further let $T$ be an $n$ by $n$
matrix of the form
$$
\beginpicture
\setcoordinatesystem units <.5cm,.3cm>
    \putrule from -.2 -4 to -.2  5
    \putrule from -.2 -4 to 0   -4
    \putrule from -.2  5 to 0    5
\putrule from 7.2 -4 to 7.2 5
\putrule from 7.2 -4 to 7 -4
\putrule from 7.2  5 to 7 5
\putrule from 0 0 to 7 0
\putrule from 3.5 0 to 3.5 5
      \put {$T'$} at 3.5 -2
      \put {$1_{n-\ell}$} at 1.5 3.2
      \put {$0$} at 5.2   3.2
      \put {$=$} at -.7     0
\putrule from -2.6 -4 to -2.6 5
\putrule from -2.6 -4 to -2.4 -4
\putrule from -2.6 5 to -2.4 5
\putrule from -1.1 -4 to -1.1 5
\putrule from -1.1 -4 to -1.3 -4
\putrule from -1.1 5 to -1.3 5
\put {=} at -3.1 0
\put {$\underline{T}_1$} at -1.8 3.2
\put {$\vdots$} at -1.8 .5
\put {$\underline{T}_n$} at -1.8 -2.5
\put { T} at -3.9 .2
\endpicture
$$
with rows $\underline{T}_i$, where $T'$ is a generic $\ell$
by $n$ matrix over $R'$.  Set $\p = T\psi$.

There exists a local ring $R''$ obtained from $R'$ by another purely
transcendental residue field extension and a prime ideal $Q$ of
$A = R''[T]$ with $\text{\rm det}(T)\notin Q$ and $\im \subset Q$ having 
the following
property: 
for $n-\ell+1\leq t\leq n$ there exist $A_Q$-regular
sequences $\underline{a}_{t}$ each of length
$n-\ell$, so that the following holds: given an arbitrary (possibly empty)
set $\Lambda = \{t_1,...,t_d\}$ of integers
$n-\ell+1\leq t_1< \dots < t_d\leq n$ and  writing $B = R''[\{\underline{T}_i|\,
i\notin \Lambda\}]_{Q^c}, \,  \underline{a} =
\underline{a}_{t_1},\dots,\underline{a}_{t_d},$ and $S = A_Q/(\underline{a})$:
\roster
\item The homomorphism $R\ra B$ is local (and regular), and the
homomorphism $B\ra S$ is local flat with a complete intersection closed
fiber.
\item $\p_S$ is a presentation matrix of $M_S$, and for $n-\ell+1\leq i\leq n$,
the $i^{\text{th}}$ row ideal $J_i$ of $\p_S$ has height at most $r-n+\ell 
= \ell-e$ if   $i\in \Lambda$ and at most $r$ otherwise.
\item $\underline{a}$ form a regular sequence on $A_Q/IA_Q$ for every 
proper ideal $I$ of $R$.
\item $\text{\rm ecodim}(S_P) = \text{\rm ecodim}(R_{P \cap R})$ for every 
prime ideal $P$ of $S$  with $P \not\in V(J_{t_1} \cdots J_{t_d})$.
\endroster
\endproclaim

\demo{Proof} Write
$U = R'v_{n-\ell+1}+ \cdots + R' v_n \subset M_{R'}$.  Since $U$ is 
generated by $\ell$ generic elements of $M$ it follows that $U$ is a 
minimal reduction of $M_{R'}$.  To simplify notation we write $R$ instead 
of $R'$ from now on.

For $n-\ell+1\leq i\leq n$ let  $\A_i$ be the ideal of
$A_i = R[\underline{T}_i]$ generated by the $i^{\text{th}}$ row ideal of $\p$.
We obtain isomorphisms
$$A_i/\A_i\cong \text{Sym}(M)$$
sending the $(i,j)$ entry $T_{ij}$ of $T$ to
$v_j$.  Since $\text{Sym}(M)$ maps onto $\R (M)$, there are $A_i$-ideals 
$\B_i$ containing $\A_i$ such that $A_i/\B_i\cong \R (M)$. Observe that
bight$(\B_i) = \text{bight}(\B_iA) = r$ (see the remarks at the start of this
section).

Let $(R'',\im'')$ be the local ring obtained from $R=R'$ by a purely
transcendental residue field extension of transcendence degree 
$(n-\ell)(\sum\mu(\B_i))$,  let $k'' = R''/\im'', ~~ \E = \otimes ^{\ell}_R 
\R (U)$ \linebreak $\otimes_R k''$ and  $\F = \otimes^{\ell}_R \R(M) 
\otimes_R k''$. Note that $\E$ is isomorphic to a polynomial ring over
$k''$ in $\ell^2$ variables which are the images of 
$T_{ij}$ in $\E$ for $n-\ell+1 \le i \le n$ and $n-\ell+1 \le j \le
n$.

  The above isomorphisms induce an isomorphism
$$A/(\im,\B_{n-\ell+1} , \dots, \B_n) \cong \F.$$
Moreover, the natural map of $k''$-algebras $\E \ra \F$ is module finite 
since $U$ is a reduction of $M$.  Its image is generated by the images in 
$\F$ of $T_{ij}$ for $n-\ell+1 \le i \le n$ and $n-\ell+1 \le j \le 
n$.  Hence these elements of $\F$ are algebraically independent over $k''$, 
because $\text{dim} \F = \ell^2$.  It follows that the image of $\Delta = 
\text{det}(T)$ in $\F$ is not nilpotent.  Thus there exists a prime ideal 
$Q$ of $A$ with $\Delta \not \in Q$ and $(\im,\B_{n-\ell+1}, \dots, \B_n) 
\subset Q$.

For every $t, n-\ell+1\leq t \leq n$, let $\underline{a}_t \subset 
R''[\underline{T}_{t}]$ be a sequence of $n-\ell$ generic elements for 
$\B_{t}\subset R[\underline{T}_{t}]$, defined using indeterminates over $R 
= R'$ as coefficients.  Such sequences exist by the definition of $R''$. As 
$(\im, \B_{t})A_{t}/\im
A_{t}$ is an ideal in a polynomial ring over a field of height
$\text{dim}(A_{t}/\im A_{t})-
\ell(M) = n-\ell$, it follows that $\underline{a}_{t}$ form a
regular sequence on $A_{t}/\im A_{t} \otimes_k k''$.

We are now ready to verify statements (1) -- (4) in the theorem. Write 
$\fin$ for the maximal ideal of $B$. As $\im \subset Q$ we have that $\im 
\subset \fin$ and thus the map $R \ra B$ is a
(regular) local homomorphism. Furthermore, $\underline{a}\subset Q$ and
$A_Q/\fin A_Q$ is flat over $(A_{t_1}/(\im
A_{t_1}))\otimes_{k} \cdots  \otimes_{k} (A_{t_d}/(\im A_{t_d})) \otimes_k 
k''$.  Thus
$\underline{a}$ form a regular sequence on $A_Q/\fin A_Q$, the closed
fiber of the map $B\ra A_Q$. Consequently, the (local) homomorphism $B\ra S =
A_Q/(\underline{a})$ is flat with complete intersection closed fiber, and
$\underline{a}$ form a regular sequence on $A_Q/IA_Q$ for any $R$-ideal $I 
\subset \im$  (\cite{Ma, p. 177}). This proves (1) and (3).

To show (2), observe that the image of $\Delta$ is a unit in $S$ since
$\Delta\notin Q$.
Thus $\p_S$ is a presentation matrix of $M_S$.  Obviously $J_i = 
\A_iS\subset \B_iS$. If $i\notin
\Lambda$ then $S$ is flat over $A_i$ and hence $\hgt(J_i)\leq \hgt(\B_iS) =
\hgt(\B_i (A_{i})_{Q\cap A_i})
\leq \bight(\B_i) = r$. If on the other hand $i\in \Delta$ then
$\underline{a}_i\subset \B_i$,
which together with the $A_Q$-regularity of $\underline{a}$ gives
$\hgt(J_i)\leq \hgt(\B_iS) \le \text{dim}(S) - \text{dim}(S/\B_i S) =
\text{dim}(A_Q) - d(n-\ell) - (\text{dim}(A_Q) - (d-1)(n-\ell) - \hgt(\B_iA_Q))
= \hgt(\B_iA_Q) - (n-\ell)\leq \bight(\B_i A) - n + \ell = r-n+\ell = \ell-e$.
This proves (2).

Finally, to show (4) notice that if $P \in \text{Spec}(S) \backslash 
V(J_{t_1} \cdot \dots \cdot J_{t_d})$, then $P \not\in V(\B_t S)$ for every 
$t \in \Lambda$.  Thus by the generic choice of $\underline{a}_t$ in 
$\B_t$, the ring $S_P$ is a localization of a polynomial ring over 
$R_{P\cap R}$.  \qed
\enddemo

\bigskip

We will often apply Theorem 2.2 in conjunction with the following 
generalization of a theorem of Serre:

\proclaim{Lemma 2.3} 
{\roster
\item"(1)" Let $f: (A,\im)\ra (B,\n)$ be a local homomorphism of
equidimensional and universally catenary Noetherian local rings, with $A$ regular.
Let $I$ be an ideal of
$A$ and $J$ be an ideal of $B$. Then
$$ \text{\rm ht}(IB + J) \leq \text{\rm ht}(I) + \text{\rm ht}(J).$$

\item"(2)" Let $B \ra S $ be a local  homomorphism of 
Noetherian local rings with $S$ equidimensional and universally catenary, 
let $K$ be an ideal of $B$, and let $J$ be an ideal of $S$.  Then 
$\text{\rm ht}(J+KS)\leq \text{\rm ht}(J) + \text{\rm ht}(K) + \text{\rm 
ecodim}(B)$.
\endroster}
\endproclaim

\demo{Proof} We first prove (1). 
Suppose first that $f$ is onto, and write $B = A/K$. Lift $J$ to
an ideal $L$ in $A$, so that $J = L/K$. Since $B$ is equidimensional
and universally catenary, $\hgt(IB + J) = \hgt((I+L)/K) =
\hgt(I+L) - \hgt(K)\leq \hgt(I) + \hgt(L) - \hgt(K) = \hgt(I) +\hgt(J)$,
where the middle inequality follows from the subaddivitity of height
in regular local rings [S, Chap. V, Thm. 3].

We now do the general case. Without loss of generality, we may assume
both $A$ and $B$ are complete: our assumptions do not change (see [Ma, 31.7]), nor
does the conclusion. We use a Cohen factorization of $f$ as in
[AFH, Thm 1.1]. There is a factorization of $f$,
$A\overset{g}\to\ra C\overset{h}\to\ra B$ where
$C$ is local, $g$ is flat and $C/\im C$ is regular, and $h$ is surjective.
Since $g$ is flat with regular closed fiber and $A$ is regular, it follows
that $C$ is also regular by [Ma, Thm. 23.7]. As $C$ maps onto $B$, to finish the proof it
suffices to prove that $\hgt(I)\geq \hgt(IC)$. However, as $A$ and $C$ are
regular and $g$ is flat,
$$\hgt(I) = \dim(A) - \dim(A/I) = \dim(C) - \dim(C/\im C) -(\dim(C/IC) - \dim(C/\im C))
= \hgt(IC).$$

We prove (2). We can pass to the completions of $B$ and $S$ and assume both rings
are complete.
Write $B = A/I$, where $A$ is a regular local ring, and lift $K$ to an
ideal $L$ in $B$, so that $K = L/I$. Note that $LS = KS$, so to prove 
(2), it is enough to prove that $\hgt(K) + \text{ecodim}(B) = \hgt(L)$
and then apply (1). But this equality is immediate.\qed
\enddemo

\bigskip

Next we give a short proof of a modified version of Theorem 2.2.  It 
requires the following definition:

\remark{Definition 2.4} Let $R$ be a Noetherian local ring with residue field $k$
(or a positively graded $k$-algebra), let $M$ be a finitely generated (graded)
$R$-module having a rank, and write $\Cal R = \Cal R(M)$. We set
$$s(M) = \text{dim}_k[(\Cal R\otimes_Rk)/\sqrt 0]_1.$$
\endremark

\bigskip

\remark{Remark 2.5} Observe that in general $\ell(M)\leq s(M)\leq \mu(M)$. If $M$
is graded and generated by forms of the same degree, then $\Cal R\otimes_Rk$
embeds into a polynomial ring over $R$ and therefore $s(M) = \mu(M)$ as
long as
$R$ is reduced and $M$ is torsionfree.
\endremark

\bigskip

\proclaim{Theorem 2.6} Let $R$ be an equidimensional universally catenary
Noetherian local ring
with algebraically closed residue field $k$, let $M$ be a finitely generated
$R$-module with rank $e$, and write $r = \mu(M) - e, ~~ s = s(M)$.
There exists a minimal presentation matrix of $M$ that has $s$ row ideals
of height at most $r$.
\endproclaim

\demo{Proof} Write $\R$ for the Rees algebra of $M$ and set
$V = [(\R\otimes_Rk)/\sqrt 0]_1$, which we identify with
affine space of dimension $s$. Consider the closed subset $X$ of $V$
whose coordinate ring is the homogeneous $k$-algebra $(\R\otimes_Rk)/\sqrt
0$.
Since $k$ is algebraically closed there exists a basis $v_1,\dots,v_s$ of $V$
contained in $X$, and then the lines $kv_1,\dots,kv_s$ all lie on $X$.

Let $z_1,\dots,z_n$ be a minimal generating set of $M$ chosen so that
$z_i$ maps to $v_i$ for $1\leq i\leq s$, and let $\p$ be a presentation
matrix with respect to $z_1,\dots,z_n$. Set $J_i$ equal to the ideal
generated by the $i^{\text{th}}$ row of $\p$. We claim $\hgt(J_i)\leq r$ for
$1\leq i\leq s$.

Let $A = R[T_1,\dots,T_n]$ be a polynomial ring, let $\im$ denote the 
maximal ideal of $R$, and for $1 \le i \le s$ consider the prime ideals $Q_i =
(\im, T_1,\dots,\hat{T_i},\dots,T_n)$ of $A$. Mapping
$T_j$ to $z_j$ for $1\leq j\leq n$,  we obtain presentations
Sym$(M)\cong A/\A$ and $\R\cong A/\B$, where $\A\subset \B$ are $A$-ideals.
As $X$ contains the line $kv_i$, we have $\B\subset Q_i$ for $1\leq i\leq s$.
Thus $\hgt(\A_{Q_i})\leq \hgt(\B_{Q_i})\leq \bight(\B) = r$. Let
$\pi_i: A_{Q_i}\lra R(T_i)$ be the $R[T_i]$-epimorphism whose kernel is
generated by the $A_{Q_i}$-regular sequence $T_1,\dots,\hat{T_i},\dots,T_n$.
Since $\hgt(\pi_i(\A_{Q_i})) + n - 1 = \hgt(\pi_i(\A_{Q_i}), T_1,...,\hat{T_i},\dots,T_n)
= \hgt(\A_{Q_i}, T_1,...,\hat{T_i},\dots,T_n)\leq \hgt(\A_{Q_i}) + n - 1$, it follows
that $\hgt(\pi_i(\A_{Q_i}))\leq \hgt(\A_{Q_i})$. But $\pi_i(\A_{Q_i}) =
J_iR(T_i)$, which gives $\hgt(J_i)\leq r$.
\qed\enddemo

\bigskip

We finish the section with two immediate consequences of Theorem 2.6.  Both 
are first height estimates for ideals of minors of matrices, stated more 
conveniently in terms of Fitting ideals of modules.

\proclaim{Corollary 2.7} Let $R$  be a regular local ring with perfect
residue field $k$ and let $M$ be a finitely generated $R$-module of rank 
$e$, and write $r = \mu(M) - e,~~ s = s(M)$. For every $1\leq i\leq s$,
$$\text{\rm ht} (\text{\rm Fitt}_{i-1}(M))\leq ir.$$
\endproclaim

\demo{Proof} There exists a flat local homomorphism $R\ra S$ where
$S$ is a regular local ring with algebraically closed residue field $K$
\cite{G,(10.3)}. Since $S$ is flat over
$R$ and $k$ is perfect, one has that
$(\Cal R(M)\otimes_SK)/\sqrt 0\cong (\Cal R(M)\otimes_Rk)/\sqrt 0)\otimes_kK$
and therefore $s(M) = s(M_S)$. We replace
$R$ and $M$ by $S$ and $M_S$, and assume that $k$ is algebraically closed.

By Theorem 2.6 there exists a minimal presentation matrix of $M$ that
has $i$ row ideals $J_1,\dots,J_i$ of height at most $r$. As
$\text{Fitt}_{i-1}(M)\subset
J_1 + \cdots + J_i$ and $R$ is a regular local ring, we conclude that
$\hgt (\text{Fitt}_{i-1}(M))\leq ir$. \qed\enddemo

\bigskip

\proclaim{Corollary 2.8} Let $R$ be a polynomial ring over a field, let
$M$ be a torsionfree graded $R$-module of rank $e$ minimally
generated by $n$ homogeneous elements of the same degree, and
write $r = n-e$. For every $1\leq i\leq n$, $\text{\rm ht} (\text{\rm 
Fitt}_{i-1}(M))\leq ir$.
In particular, for every submodule $U$ of $M$ generated by $t < n$
elements, $\text{\rm ht} (U:_RM) \leq (t+1)r$.
\endproclaim

\demo{Proof} We may assume that the ground field is perfect. Writing
$\im$ for the irrelevant maximal ideal of $R$ we observe that $s(M_{\im})= 
s(M) = n$
and $(U:_RM)_{\im}\subset \text{Fitt}_t(M_{\im})$. The assertions now
follow from
Corollary 2.7. \qed
\enddemo

\bigskip

\bigskip
\bigskip
\centerline{\bf 3. Heights of Determinantal Ideals\rm}
\bigskip

The classical theorem of Bruns ([B1, Cor. 1]) states that in a Noetherian 
ring $R$, the height of the ideal of $i$ by $i$ minors of an $n$ by $m$ 
matrix of rank $r$ cannot exceed the ``generic'' value $N(i,r,m,n)$ defined 
as follows:  let $X$ be a generic $n$ by $m$ matrix and set $N(i,r,m,n) := 
\hgt(I_i(X)) - \hgt(I_{r+1}(X)) = (r-i+1)(m+n-r-i+1)$.  This is exactly the 
height of the  ideal of $i$ by $i$ minors of the image of $X$ in the ring 
$R[X]/I_{r+1}(X)$ (note the image of $X$ has rank $r$ in this 
ring).  However, if we also insist that the base ring $R$ be regular, then 
it is by no means clear that this maximum is ever attained.  The main 
results known for the regular case are due to Bruns and Faltings ([B1, Thm. 
3], [F]), and their results apply only to the case $i=r$.  In Corollary 
3.6.1 below we establish a bound for the height of the ideal of $i$ by $i$ 
minors of an $n$ by $m$ matrix of rank $r$ over a regular ring that is 
roughly $(r-i)(\max\{m,n\}-i+1) + \max\{m-i+1,n\}$.

A second, related problem is to estimate the height of the ideal of $i$ by 
$i$ minors modulo the ideal of $i+1$ by $i+1$ minors.  Again, the best 
general bound is $N(i,i,m,n) = m+n-2i+1$, but one may expect better results 
if $R$ is regular and the rank $r$ of the matrix is not maximal.  We 
address this issue in Corollary 3.9.1, where the bound $\max\{m-i+1,n\} +r 
-i$ is established.

Both problems are special cases of the following, more general 
question:  How can one estimate the height of the ideal of $i$ by $i$ 
minors of a matrix of rank $\bar{r}$ that can be ``lifted'' to a matrix of 
rank $r$ over a ring $R$?
Theorem 3.1, the main result of this section, gives such a bound involving 
the difference $r - \bar{r}$ of the ranks and the embedding codimension of 
$R$.  The proof of this result relies on the work of Section 2 about row 
ideals of small height. The theorem gives particularly strong estimates if 
the matrix can be lifted in such a way that the increase in the rank is 
compensated by a decrease in the embedding codimension of the ambient ring.

\bigskip

\proclaim{Theorem 3.1} Let $R$ be an equidimensional universally catenary 
Noetherian local ring, let $\var$ be an $n$ by $m$ matrix of rank $r$
with entries in $R$, and let $I$ be an $R$-ideal.  Assume that  $M 
= \text{\rm Coker}(\var)$ has a rank, and write $\ell = \ell(M), ~ \bar{R} 
= R/I, ~ \bar{\var} = \var_{\bar{R}}, ~ \bar{r} = \text{\rm 
rank}(\bar{\var})$.  Let $i \le \bar{r}$ be an integer so that 
$I_i(\bar{\var}) \neq \bar{R}$.  Set $\d = r - \bar{r}$ and $\e = 
\max_P\{\text{\rm ecodim}(R_P)\}$, where the maximum is taken over all 
prime ideals $P$ of $R$ not containing $I_i(\var)$.
\roster
\item
$$
  \text{\rm ht}(I_i(\bar{\var})) \le 
\max\{(\min\{n-\ell,\bar{r}\}-i+1)(m-i+1 +\max\{0, n-\ell-\bar{r}\}),$$
$$\qquad\qquad(\bar{r}-i)(\max\{m,n+\e\} - i + \d + 1) + \ell + \d + 
\text{\rm ecodim}(R)\}$$
$$\quad \le  (\bar{r}-i)(\max\{m,n+\e\} - i + \d + 1) + \max\{m-i+1, \ell + 
\text{\rm ecodim}(R)\}+ \d.$$

\medskip

\item  If the $\bar{R}$-module $\bar{M} = M_{\bar{R}}$ is not a direct sum 
of a torsion module and a free module, $\bar{M}_{\bar{P}}$ is free for 
every prime $\bar{P}$ of $\bar{R}$ with depth$(\bar{R}_{\bar{P}}) \le 1$ 
and $M_P$ is of linear type for every associated prime $P$ of $I$, then
$$\text{\rm ht}(I_i(\bar{\var})) \le (\bar{r}-i)(\max\{m,n + \e\} - i + \d 
+ 1) + \ell + \d + \text{\rm ecodim}(R).$$
\endroster
\endproclaim

\bigskip

Before proving the theorem we wish to make several comments. First notice 
that $\e = 0$ in case $R$ is locally regular on the punctured spectrum.  If 
the $\bar{R}$-module $\bar{M}$ is a direct sum of a torsion module and a 
free module then trivially $\hgt(I_i(\bar{\var})) \le (\bar{r}-i+1) 
(m-i+1)$.  It is also obvious that one can replace the bound of part (1) by 
the better formula of (2) whenever $i \ge n - \ell + 1$.  Finally, the 
estimates of Theorem 3.1 are sharp for $\var$ a generic matrix with entries 
in the localization of a polynomial ring over a regular ring and $I = 
I_{\bar{r}+1} (\var)$, if $n \le m$ or $i = 1$.

\bigskip

\demo{Proof of Theorem 3.1} We first prove that the second inequality of
(1) is true, namely that
$$
\max\{(\min\{n-\ell,\bar{r}\}-i+1)(m-i+1 +\max\{0, n-\ell-\bar{r}\}),$$
$$\qquad\qquad(\bar{r}-i)(\max\{m,n+\e\} - i + \d + 1) + \ell + \d +
\text{\rm ecodim}(R)\}$$
$$\quad \le  (\bar{r}-i)(\max\{m,n+\e\} - i + \d + 1) + \max\{m-i+1, \ell +
\text{\rm ecodim}(R)\}+ \d.$$

We prove each term in the maximum on the left hand
side of the inequality is at most the right hand side. This is clear for the
second term. It remains to see why
$$(\min\{n-\ell,\bar{r}\}-i+1)(m-i+1 +\max\{0, n-\ell-\bar{r}\})$$
$$\leq  (\bar{r}-i)(\max\{m,n+\e\} - i + \d + 1) + \max\{m-i+1, \ell +
\text{\rm ecodim}(R)\}+ \d.$$ 
By possibly lessening the right hand side and increasing the left-hand side,
it is enough to prove that
$$ (\bar{r}-i+1)(m-i+1 +\max\{0, n-\ell-\bar{r}\})\leq 
(\bar{r}-i)(m-i+\d+1) + m-i+1+\d = (\bar{r}-i+1)(m-i+\d+1),$$
and for this it suffices to prove that $\max\{0, n-\ell-\bar{r}\}\leq
\d = r-\bar{r}$. Clearly $0\leq \d$. The inequality
$n-\ell-\bar{r}\leq r-\bar{r}$ is equivalent to the inequality
$n-r\leq \ell$, which is always true, since $n - r = e = \text{rank}(M)\leq \ell$.

We use induction on $n$ to prove the first inequality of Theorem 3.1.
Suppose that $n = 1$. In this case, $M = R/J$, where $J$ is an ideal
with $m$-generators. By assumption, $M$ has a rank, which is of necessity
either $0$ or $1$. However the rank cannot be $1$, since then $J_P = 0$ for
all associated primes of $R$, and hence $J = 0$ and $M = R$ is free.
 Thus the rank of $M$ is $0$, and then $J$ contains a
non-zerodivisor. It follows that the analytic spread of $M$ is $0$, since
we always mod out torsion to compute the analytic spread. Hence,
$\ell = 0, r = 1, n=1$,  and $\bar{r}$ is either $0$ or $1$. 
If $\bar{r} = 0$, then the theorem is vacuous. Hence we may assume that
$\bar{r} = 1$ also, and $i = 1$. In this case the inequality reads:
$$\text{\rm ht}(I_1(\bar{\var})) \le \max\{m, \text{ecodim}(R)\}.$$
By the Krull height theorem, the height of $I_1(\bar{\var})$ is at most its
number of generators, which is bounded by $m$, proving the case $n = 1$.

  We may assume that the entries of $\var$ lie 
in the maximal ideal of $R$.  We claim that we may assume that
$I = P$ is a prime ideal. Let $P$ be  
a minimal prime of $I$ having maximal dimension. We write $r_P$ for the
rank of $\phi_{R/P}$. There are three cases,
depending on the relationship of $r_P$ to $i$ and $\bar{r}$. Note that
$r_P\leq \bar{r}$.

Case 1. $r_P = \bar{r}$. Since $R$ is equidimensional and
catenary, $\text{\rm ht}(I_i(\bar{\var})) = \text{\rm ht}(I_i(\var_{R/P}))$.
Hence the left-hand side of in the inequality of (1) doesn't change, but
neither does the right-hand side in this case.

Case 2. $r_P < i$. Then $I_i(\var)\subseteq P$, and 
$\text{\rm ht}(I_i(\var_{R/P})) = 0$. Since the right-hand side of (1)
is nonnegative, the inequality holds.

Case 3. $i\leq r_P$. In this case we prove that as a function of
$\bar{r}$,  the right-hand side
of (1) is nonincreasing as we decrease 
$\bar{r}$ to  $i$. Since $i\leq r_P\leq \bar{r}$ and since 
 $\text{\rm ht}(I_i(\bar{\var})) = \text{\rm ht}(I_i(\var_{R/P}))$,
this will prove our claim. The right-hand side of (1) is a maximum of
two terms. Decreasing $\bar{R}$ by one changes the second term,
$(\bar{r}-i)(\max\{m,n+\e\} - i + \d + 1) + \ell + \d +
\text{\rm ecodim}(R)$,  to
$(\bar{r}-i-1)(\max\{m,n+\e\} - i + \d + 2) + \ell + \d + 1 +
\text{\rm ecodim}(R)$. Subtracting the first from the second gives
the value max$\{m,n+\e\} + r + 1 - 2\bar{r}$, which is always
nonnegative.

The first term, $(\min\{n-\ell,\bar{r}\}-i+1)(m-i+1 +\max\{0, n-\ell-\bar{r}\})$
can only increase if $n-\ell-\bar{r}\geq 0$. Then as $\bar{r}$ decreases by
$1$, $\max\{0, n-\ell-\bar{r}\}$ will increase by $1$. However,
in this case $\min\{n-\ell, \bar{r}\}$ will be $\bar{r}$ and will decrease
by $1$. Then the product has the form 
$(\bar{r}-i+1)(m-i+1 + (n-\ell-\bar{r}))$, and when we replace
$\bar{r}$ by $\bar{r}-1$ we obtain
$(\bar{r}-i)(m-i+1 + (n-\ell-\bar{r}) +1)$. But 
$(\bar{r}-i+1)(m-i+1 + (n-\ell-\bar{r}))\geq (\bar{r}-i)(m-i+1 + (n-\ell-\bar{r}) +1)$
since $2\bar{r}\leq m + n -\ell +1$. 

Thus
 we may suppose that 
$\bar{R}$ is a domain, hence equidimensional. 
We use the notation of Theorem 2.2 and in addition set $\underline{a}_j = 
0$ whenever $j \le n - \ell$.  For $0 \le j \le n$ let $\p_j$ be the $j$ by 
$m$ matrix consisting of the first $j$ rows of $\phi$, and define
$$t = \min\{j\mid I_i(\phi) \subset \sqrt{(I_i(\phi_j),I, 
\underline{a}_j)_Q}\}.$$
We may assume that $i \le t$. For suppose that $t < i$. Then
$I_i(\phi_t) = 0$ so we would have that
$I_i(\phi) \subset \sqrt{(I,\underline{a}_j)_Q}$. The map from $R$ to
$S$ is flat, and the $\underline{a}_j$ from a regular sequence in
$A_Q/IA_Q$. Hence if $s\in I_i(\phi)$, then for large $N$,
 $s^N\in (I,\underline{a}_j)_Q\cap R = I$,
the last equality by flatness. Since $I$ is
prime, we obtain that $s\in I$, and then $I_i(\phi)\subseteq IR$ and we
are done. Henceforth we assume that $i\le t$. 

  We apply Theorem 2.2 with $\L = \emptyset$ if $t \le 
n - \ell$ and $\L  = \{t\}$ if $t \ge n - \ell + 1$.  Let $J_t$ be the 
$t^{\text{th}}$ row ideal of the matrix $\phi_S$, and write $\bar{S} = 
S/IS, ~ \bar{J}_t = J_t \bar{S}$.  By Theorem 2.2, $R \subset S$ and 
$\bar{R} \subset \bar{S}$ are flat local extensions, $S$ and $\bar{S}$ are 
equidimensional and catenary, and ecodim$(S_P) \le \e$ for every prime $P$ 
of $S$ not containing $I_i(\var) \cdot J_t$.  Notice that 
$I_i(\phi_{\bar{S}}) \subset \sqrt{I_i((\phi_t)_{\bar{S}})}$ and $I_i 
(\phi_{\bar{S}}) \not\subset \sqrt{I_i((\phi_{t-1})_{\bar{S}})}$ according 
to the definition of $t$.  Again by Theorem 2.2, $\hgt(J_t) \le r-n+\ell$ 
if $t \ge n - \ell+1$.  Furthermore as $I_i((\phi_{t-1})_S) + IS$ is 
extended from $B$, Lemma 2.3 implies that
$$ \hgt
(J_t + I_i((\phi_{t-1})_S) + IS) \le \hgt (J_t) + \hgt(I_i((\phi_{t-1})_S) 
+ IS) + \text{ecodim}(R). $$
Thus by our equidimensionality conditions,
$$\hgt(\bar{J}_t + I_i((\phi_{t-1})_{\bar{S}})) \le \hgt(J_t) + 
\hgt(I_i((\phi_{t-1})_{\bar{S}})) + \text{ecodim}(R).$$
Since $\bar{S}$ is flat over $\bar{R}$ and
$$ I_i(\var_{\bar{S}}) = I_i(\phi_{\bar{S}}) \subset 
\sqrt{I_i((\phi_t)_{\bar{S}})} \subset \sqrt{\bar{J}_t + I_i(( 
\phi_{t-1})_{\bar{S}})},$$
we conclude that
$$\hgt(I_i(\bar{\var})) = \hgt(I_i(\var_{\bar{S}})) \le \hgt(J_t) + 
\hgt(I_i(( \phi_{t-1})_{\bar{S}})) + \text{ecodim}(R).$$
To simplify notation we will henceforth write $\phi, \phi_j, \bar{\phi}, 
\bar{\phi}_j$ instead of $\phi_S, (\phi_j)_S, \phi_{\bar{S}}, 
(\phi_j)_{\bar{S}}$.  With this we have
$$
\sqrt{I_i(\bar{\phi}_{t-1})} \subsetneq \sqrt{I_i(\bar{\phi}_t)}\leqno(3.2)$$
and
$$ \hgt(I_i(\bar{\var})) = \hgt(I_i(\bar{\phi}_t)) \le \hgt(J_t) + 
\hgt(I_i(\bar{\phi}_{t-1})) + \text{ecodim}(R). \leqno(3.3) $$

\bigskip

\medskip

Case 1:  $t \le n - \ell$.  In this case $I_i(\var_{\overline{S}}) =
I_i(\phi_{\overline{S}})\subseteq \sqrt{I_i(\bar{\phi}_t)_{\overline{S}}}$, so
that  $\hgt(I_i(\bar{\var})) \le \hgt (I_i(\phi_{n- \ell}))
\le (n-\ell - i + 1)(m-i+1)$,
and according to [B1, Cor. 1], 
$$ \hgt(I_i(\bar{\var})) \le (\bar{r}-i+1)(m+n-\ell-\bar{r} - i + 1)$$
The first inequality of (1) follows, and the second holds because $\ell 
=\ell(M) \ge \text{rank}(M) = n-r$.

\bigskip

Case 2:  $t \ge n - \ell + 1$.  In this case $\hgt (J_t) \le r - n + \ell$, 
and therefore (3.3) yields
$$ \hgt(I_i(\bar{\var})) \le r - n + \ell + \hgt(I_i(\bar{\phi}_{t-1})) + 
\text{ecodim}(R).
\leqno(3.4) $$
By (3.2) there exists a prime ideal $P$ of $S$ with $I_i(\phi_{t-1}) + IS 
\subset P$ and $I_i(\phi_t) \not\subset P$.  Since $I_i(\phi_t)$ is 
contained in $I_{i-1}(\phi_{t-1})$, in $I_i(\var)S$, and in $J_t + 
I_i(\phi_{t-1})$, one automatically has $I_{i-1}(\phi_{t-1}) \not\subset P$ 
as well as $I_i(\var) \cdot J_t \not\subset P$.  By the latter, 
$\text{ecodim}(S_P) \le \e$.  Set $s = \max\{j | I_j(\phi) \not\subset 
P\}$.  Clearly $1 \le i \le s$.  Recall that $I_{i-1} (\phi_{t-1})_P = S_P$ 
and $I_i(\bar{\phi}_{t-1})_P \neq \bar{S}_P$. Thus without changing the 
ideal $I_i(\bar{\phi}_{t-1})_P$, we may perform elementary row and column 
operations over $S_P$ to assume that

\bigskip

$$\beginpicture
\setcoordinatesystem units <.5cm,.4cm>
\putrule from -.2 -5 to -.2 5
\putrule from -.2 -5 to 0 -5
\putrule from -.2 5 to 0 5
\putrule from 2.5 5 to 2.5 -5
\putrule from 0 2 to 7.5 2
\putrule from 2.5 -1.5 to 7.5 -1.5
\putrule from 5 2 to 5 -5
\putrule from 7.7 -5 to 7.7 5
\putrule from 7.7 -5 to 7.5 -5
\putrule from 7.7 5 to 7.5 5
\put {=} at -.7 0
\put {$1_{i-1}$} at 1.2 3.5
\put {$0$} at 5 3.5
\put {$0$} at 1.3 -1.5
\put {$0$} at 3.8 -3.5
\put {$\p'$} at 3.8 0
\put {$1_{s-i+1}$} at 6.2 -3.5
\put {$\p''$} at 6.3 0
\put {$\p_P$} at -1.7 .1
\endpicture
$$

\bigskip
\bigskip

\noindent where $\phi',\phi''$ have entries in the maximal ideal of 
$S_P$.  Notice that the $n-s$ by $m-s$ matrix $\phi'$ has rank $r-s$ and 
$\bar{\phi'}$ has rank $\bar{r} - s$, with $\bar{\phi'}, \bar{\phi''}$ 
standing for $\phi'_{\bar{S}_P},  \phi''_{\bar{S}_P}.$

Since $I_i(\bar{\phi}_{t-1})_P \subset I_1(\bar{\phi'}) + I_1(\bar{\phi''}) 
\neq \bar{S}_P$ and $\bar{S}_P$ is equidimensional, we obtain
$$ \hgt(I_i(\bar{\phi}_{t-1})) \le \hgt (I_i(\bar{\phi}_{t-1})_P) \le \mu 
(I_1(\bar{\phi}'')) + \hgt(I_1(\bar{\phi}'))$$
$$ \qquad \le (s - i + 1)(n-s) + \hgt(I_1(\bar{\phi}')).$$
Thus by (3.4),
$$ \hgt (I_i(\bar{\var})) \le r - n + \ell + (s-i+1)(n-s) + 
\hgt(I_1(\bar{\phi}')) + \text{ecodim}(R).
\leqno(3.5) $$

\bigskip

Applying the induction hypothesis to the matrix $\phi'$ yields

%\pagebreak

$$~~~~~\hgt(I_1(\phi'))   \le ((\bar{r}-s)-1)(\max\{m-s,(n-s) + \e\} - 1 + 
\d + 1)$$
$$  + \max\{(m-s) - 1 + 1,(n-s) + \e\} + \d $$
$$  = (\bar{r}-s)(\max\{m,n+\e\} - s + \d).$$
Hence by (3.5),
$$ \hgt(I_i(\bar{\var})) \le r - n + \ell + (s-i+1)(n-s) + 
(\bar{r}-s)(\max\{m,n + \e\} - s + \d) + \text{ecodim}(R)$$
$$ \!\!\!\!\! \le r - n + \ell + n-i + (\bar{r}-i)(\max\{m,n + \e\} - i + 
\d) + \text{ecodim}(R),$$
because $i \le s$ and $n-s \le \max\{m,n + \e\} - i + \d$.  It follows that
$$ \hgt (I_i(\bar{\var})) \le (\bar{r} - i)(\max \{m,n + \e\} - i + \d + 1) 
+ \ell + \d + \text{ecodim}(R),$$
proving (1) in Case 2 as well.

\medskip

To show part (2) first notice that the $\bar{R}$-module $\bar{M}$ has a 
rank, as can be seen from the Abhyankar-Hartshorne connectedness lemma (see 
Hartshorne [Ha2]).  The natural map $\text{Sym}(M) \to \text{Sym}(\bar{M})$ 
induces an epimorphism $\R(M) \to \R(\bar{M})$ since $M$ is of linear type 
locally at every associated prime of $I$.  Therefore $\ell(M) \ge 
\ell(\bar{M})$.  On the other hand $\ell(\bar{M}) \ge \text{rank}(\bar{M}) 
+ 1$ by Proposition 1.3.   Therefore $\ell \ge n - \bar{r} + 1,$ and (2) 
follows from (1).
\qed\enddemo

\bigskip

\proclaim{Corollary 3.6} Let $R$ be an equidimensional universally catenary 
Noetherian local ring and  let  $\var$ be an $n$ by $m$ matrix of rank $r$ 
with entries in $R$.  Assume that $M = \text{\rm Coker}(\var)$ has a rank 
and write $\ell = \ell(M)$.  Let $i \le r$ be an integer such that 
$I_i(\var) \neq R$.  Set $\e = \max_P\{\text{\rm ecodim}(R_P)\}$ where the 
maximum is taken over all prime ideals $P$ of $R$ not containing $I_i(\var).$
\roster
\item
$$
  \text{\rm ht}(I_i(\var)) \le \max\{(n-\ell\!-\!i\!+\!1)(m\!-\!i\!+\!1), 
(r-i)(\max\{m,n+\e\}-i+1) + \ell + \text{\rm ecodim}(R)\}$$
$$ \le (r-i)(\max\{m,n+\e\} - i+1) + \max\{m-i+1, \ell +  \text{\rm 
ecodim}(R)\}.$$
\item  If  $M$ is not a direct sum of a torsion module and a free module then
$$\text{\rm ht}(I_i(\var)) \le (r-i)(\max\{m,n + \e\} - i + 1) + \ell 
+  \text{\rm ecodim}(R).$$
\endroster
\endproclaim

\demo{Proof}  Apply Theorem 3.1 with $I = 0$ and use that $\ell \ge 
\rank(M)$.
\qed\enddemo

\bigskip

In the setting of Corollary 3.6, part (1) could also be deduced from (2): 
for if $M$ is a direct sum of a torsion module and a free module then 
obviously $\hgt(I_i(\var)) \le (r-i+1)(m-i+1).$

\bigskip

\proclaim{Corollary 3.7} Let $R$ be an equidimensional universally catenary 
Noetherian local ring and  let  $\var$ be an $n$ by $m$ matrix of 
$\text{\rm rank} r$ with entries in $R$.  Assume that $M = \text{\rm 
Coker}(\var)$ has a rank and write $\ell = \ell(M).$
\roster
\item {\rm ([B1, Cor. 1])} If $M$ is not free then $\text{\rm 
ht}(I_r(\var)) \le \max\{m - r + 1, \ell + \text{\rm ecodim}(R)\}.$

\item  {\rm ([F, Kor. 1])}  If  $M$ is not a direct sum of a torsion module 
and a free module then
$\text{\rm ht}(I_r(\var)) \le  \ell +  \text{\rm ecodim}(R).$
\endroster
\endproclaim

\demo{Proof} Set $i = r$ in Corollary 3.6. \qed \enddemo

\proclaim{Corollary 3.8} Let $R$ be an equidimensional universally catenary 
Noetherian local ring  of dimension $d$ and  let  $M$ be a finitely 
generated $R$-module having a rank.  Let $\L$ be the set of all prime 
ideals $Q$ such that the $R_Q$-module $M_Q$ is not a direct sum of a 
torsion module and a free module.  If $\L$ is nonempty then
$$ d \le \max_{Q \in \L} \{\mu_Q (M) + \text{\rm ecodim} (R_Q) + \dim 
(R/Q)\}.$$
\endproclaim

\demo{Proof}  We may factor out the torsion of $M$ to assume that $M$ is 
torsionfree.  Notice this does not change the set $\L$.  Choose $Q$ minimal 
in $\L$.  Then $M_P$ is free for all primes $P \subsetneq Q$.  If $\var$ is 
a matrix minimally presenting $M_Q$ we let $r$ be the rank of $\var$.  Our 
choice of $Q$ shows that $\sqrt{I_r(\var)} = QR_Q$.  Corollary 3.7.2 then 
gives $\hgt(I_r(\var)R_Q) \le \mu_Q(M) + \text{ecodim}(R_Q)$.  Hence $d - 
\dim(R/Q) = \dim(R_Q)= \hgt(I_r(\var)R_Q) \le \mu_Q(M) + \text{ecodim}(R_Q)$, 
from which the corollary follows.
\qed\enddemo

\proclaim{Corollary 3.9} Let $R$ be an equidimensional universally catenary 
Noetherian local ring and let $\var$ be an $n$ by $m$ matrix of rank $r$ 
with entries in $R$.  Assume that $M = \text{\rm Coker}(\var)$ has a rank 
and write $\ell = \ell(M)$.  Let $i \le r$ be an integer such that 
$I_i(\var) \neq R$.
\roster
\item $\text{\rm ht}(I_i(\var)/I_{i+1}(\var)) \le \max\{m - i + 1, \ell + 
\text{\rm ecodim}(R)\} + r-i.$

\item   If  $i \ge n-\ell + 1$, then
$$\text{\rm ht}(I_i(\var)/I_{i+1}(\var)) \le  \ell +  r-i + \text{\rm 
ecodim}(R)$$
and in particular
$$\text{\rm ht}(I_i(\var)) \le (r-i+1)(\ell + \text{\rm ecodim}(R)) + 
\pmatrix
r-i+1 \\ 2
\endpmatrix.$$
\endroster
\endproclaim

\demo{Proof} Apply Theorem 3.1 with $I = I_{i+1}(\var).$  Notice that 
$\bar{r} = i$ and $n - \ell \le r$. Iterate to get the second statement.
\qed \enddemo

\bigskip

The reader may want to compare Corollary 3.9.2 to Corollary 2.7.  The 
significance of both formulas is that they do not involve $m$.

\bigskip

The above result leads to improved height bounds for $I_i(\var)$ if one 
knows a priori that for some $j \ge i$, the height of $I_j(\var)$ is 
``smaller than expected''.  Applying this observation to ideals one obtains:

\proclaim{Corollary 3.10}  Let $R$ be an equidimensional universally 
catenary Noetherian local ring with residue field $k$  and let $J$ be an 
$R$-ideal with $\text{\rm grade}(J) > 0$.  Write $g = \text{\rm ht}(J), ~ 
\ell = \ell(J), ~ n = \mu(J)$, and $m = \text{\rm dim}_k \text{\rm Tor}^R_1 
(k,J).$  Let  $i$ be an integer with $g-1 \le i \le n-1$.
\roster
\item  If $i \le \ell-1$ then
$$\text{\rm ht} (\text{\rm Fitt}_i(J)) \le (i-g+1)(\ell+g-1 +  \text{\rm 
ecodim}(R)) +
\pmatrix  i - g+1 \\2 \endpmatrix + g.$$

\item   If  $i \ge \ell$, then
$$\text{\rm ht}(\text{\rm Fitt}_i(J)) \le (\ell - g)(\ell + g-1 
+  \text{\rm ecodim}(R)) +
\pmatrix \ell - g  \\ 2 \endpmatrix  + g$$
$$ \qquad \qquad \qquad\quad + (i - \ell + 1) \max\{m - n + \ell + i , 
\frac{3 \ell + i}{2} - 1 +  \text{\rm ecodim}(R)\}.$$
\endroster
\endproclaim

\demo{Proof}  Notice that $\text{\rm ht}(\text{\rm Fitt}_{g-1}(I)) \le g$ 
and apply Corollary 3.9.
\qed \enddemo

\bigskip

\proclaim{Example 3.11}  Let $R$ be a regular local ring and let $J$ be a 
proper $R$-ideal with
$g = \text{\rm ht}(J)$ and $\ell = \ell(J)$.
\roster
\item \text{\rm (}Non-complete-intersection locus, \text{\rm [Hu, Thm. 1.1])}  If 
$J$ is not a complete intersection then
$\text{\rm ht}(\text{\rm Fitt}_g(J)) \le \ell + 2g - 1$.
\item \text{\rm (}Non-almost-complete-intersection locus\text{\rm)}  If 
$\text{\rm Ext}^g_R(J,R) = 0, ~ J_Q$ is a complete intersection for every 
prime $Q$ containing $J$ with $\text{\rm dim}(R_Q) = g$, and $J$ is not an 
almost complete intersection, then  $\text{\rm ht}(\text{\rm 
Fitt}_{g+1}(J)) \le 2 \ell + 3g - 1$.
\endroster
\endproclaim

\demo{Proof}  In (1) we may suppose that $\text{\rm ht}(\text{\rm 
Fitt}_g(J)) \ge g+1$.  But then $J$ satisfies $G_{g+1}$, and hence $\ell 
\ge g+1$ by [CN].  The assertion follows from Corollary 
3.10.1.  Likewise in (2) one can assume that $\text{\rm ht}(\text{\rm 
Fitt}_{g+1}(J)) \ge g+2$.  Thus $J$ satisfies $G_{g+2}$, and  therefore 
$\ell \ge g + 2$ according to [CEU, 4.4 and 3.4(a)] and Proposition 
1.2.  We may apply Corollary 3.10.1.
\qed \enddemo

\bigskip
\bigskip
\centerline{\bf 4.  A Family of Examples\rm}
\bigskip 

We present a class of $n$ by $m$ matrices of rank $r$ which
show that the inequalities of Corollary 3.6.2 are fairly sharp for all
values of $i, r,m,n$.  Unlike the examples given in the introduction, these
matrices have no generalized zeros.

\bigskip 

\noindent{\bf Example 4.1.}  Let $i,r,m,n$ be integers with $1 \le i \le r
\le n \le m$ and let $\var $ be the product of a generic $n$ by $r$ matrix
with a generic $r$ by $m$ matrix.  One has
$$ \hgt(I_i(\var)) = \cases (r\!-\!i\!+\!1)(n\!-\!i\!+\!1) & \text{if } m
\ge n+r-i+1 \\
(r\!-\!i\!+\!1)(n\!-\!i\!+\!1) - \frac{(r+n-m-i+1)^2}{4} &   \text{if }
r+n-m-i+1 > 0 \text{ and even} \\
(r\!-\!i\!+\!1)(n\!-\!i\!+\!1) - \frac{(r+n-m-i+1)^2-1}{4} &  \text{if }
r+n-m-i+1 > 0 \text{ and odd}.
\endcases $$
 
\bigskip
 
\demo{Proof}  We may assume that the ambient ring $R$ is obtained by
adjoining the entries of the two generic matrices to a ring  $k$. The
height of $I_i(\var)$ cannot decrease when  $k$ is replaced by the residue
field of any minimal prime of $k$, and it cannot increase if we pass to the
residue field of $P \cap k$ for some minimal prime $P$ of $I_i(\var)$
having minimal height.  Thus it suffices to consider the case where $k$ is
a field, and we may even assume that $k$ is algebraically closed.
 
Let $X$ be the closed subset of $\P^{r(m+n)-1}_k = \P(\Hom_k(k^m,k^r)
\times \Hom_k(k^r,k^n))$ defined by the homogeneous ideal
$I_i(\var)$.  Notice that
$X = \{[(\a,\b)] \mid \rank(\b\a) \le i-1\}$, where $\a \in
\Hom_k(k^m,k^r)$ and $\b \in \Hom_k(k^r,k^n)$.  For $0 \le s \le r-i+1$ set
$X_s = \{[(\a,\b)] \mid \rank(\a)\le s + i-1, ~ \rank(\b) \le r-s, ~
\rank(\b\a) \le i-1\}.$   As $X$ is the union of the closed subsets $X_s$,
our formula will follow once we have shown that
$$\dim X_s = (s+i-1)(r+m-s-i+1) + (r-s) n + (i-1)s-1.$$
 
In doing so we even show that $X_s$ is irreducible and we construct an
explicit desingularization (see also Huneke and Ulrich [1987, the proof of
3.16], and Arbarello, Cornalba, Griffiths and Harris [1985, Chapter II,
Section 2]).  Let $Y$ be the flag variety Fl$(s,s+i-1;k^r) = \{(U,V) \mid U
\subset V \subset k^r\},$ where $U$ and $V$ are subspaces of dimension $s$
and $s+i-1$, respectively.  In $Y \times \P^{r(m+n)-1}_k$ consider the
closed subset $Z = \{((U,V),[(\a,\b)])\mid \text{Image}(\a) \subset
V,$  Ker$(\b) \supset U\}.$  The projections onto the first and second
factor of $Y \times \P^{r(m+n)-1}_k$ yield surjective morphisms
%%%%%%%%%%%%%%%%%%%%%%%%%%%%%%%%%%%%%%%%%%%%%
$$
\diagram
& \!\!\!\!\!\! Z  & \\
\ldTo^f(1.5,1.9) & \qquad \qquad \rdTo^g(1.5,1.9) \\
\!Y\qquad  & & \quad  X_s
\enddiagram
$$
$Y$ is irreducible of dimension $(s+i-1)(r-s-i+1) + (i-1)s$.  The
fibers of $f$ over all closed points $(U,V)$ of $Y$ are isomorphic to
$\P(\Hom_k(k^m,V) \times \P_k(k^r/ U,k^n)) \cong \P^{m(s+i-1) +
n(r-s)-1}_k$, hence are irreducible of constant dimension.  Since,
furthermore, \linebreak $Z \subset Y \times \P^{r(m+n)-1}_k$, it follows
that $Z$ is irreducible (see Eisenbud [1995, Exercise 14.3]).  One
necessarily has
$$\dim Z \!=\! \dim Y + \dim \P^{m(s+i-1) + n(r-s) - 1}_k =
(s+i-1)(r+m-s-i+1) + (r-s) n+(i-1)s-1,$$
as can be seen, for instance, from the lemma of generic flatness (see
Eisenbud [1995, 14.4]).  On the other hand, since $Z$ is irreducible and
$g$ is surjective, $X_s$ is irreducible as well.  As $\{[\a,\b)] \mid
\rank(\a) \le s+i-2$ or rank$(\b) \le r-s-1\} \cap X_s$ is empty or a
closed proper subset of $X_s$, it follows that for every closed point
$(\a,\b)$ in some dense open subsets of $X_s$, the fiber of $g$ over
$(\a,\b)$ consists of the single point $((\text{Ker}(\b),\text{Image}(\a)),
[(\a,\b]).$  Thus again by generic flatness, $\dim X_s = \dim Z$, which
proves  our assertion.
\qed \enddemo
 
\bigskip
\bigskip
\centerline{\bf 5.  Some Results on Symmetric Matrices \rm}
\bigskip

We prove the conjecture of Problem 1 in the 
extremal cases $i=1$ and $i=n-1$ if the ring is regular.
 
\bigskip
 
\proclaim{Proposition 5.1}  Let $(R,\im)$ be a regular local ring and let
$\var$ be a symmetric $n$ by $n$ matrix of rank $r$ with entries in $\im$.
 
\roster
\item $\text{\rm ht}(I_1(\var)) \le rn - \left( \matrix
r \\2 \endmatrix \right)$.
 
\item If $2$ is a non zerodivisor on $R$ and $r=n-1$, then
$$\text{\rm ht}(I_{n-1} (\var)) \le 2.$$
\endroster
\endproclaim
 
\demo{Proof}  To prove (1) we apply Theorem 2.1 to the module $M =
\text{Coker} (\var)$.  One has $\ell(M) \ge \rank (M) = n-r$.  By the
theorem there exists a local homomorphism $R \ra S$ with $S$ a localization
of a polynomial ring over $R$, and an invertible $n$ by $n$ matrix $T$ over
$S$ so that $n-r$ row ideals $J_1, \dots, J_{n-r}$ of $\Psi = T \var T^*$
have height at most $r$.  By the symmetry of $\var, ~~ \mu(I_1(\Psi)/(J_1 +
\cdots + J_{n-r})) \le \left( \matrix r+1 \\ 2 \endmatrix
\right)$.  Therefore $\hgt(I_1(\var)) = \hgt(I_1(\var)S) = \hgt(I_1(\Psi))
\le \hgt(J_1 + \cdots + J_{n-r}) + \left( \matrix r+1 \\ 2 \endmatrix
\right) \le (n-r) r +
\left( \matrix r+1 \\ 2 \endmatrix \right) = rn -
\left( \matrix r \\ 2 \endmatrix \right).
$
 
\medskip
 
To prove (2) we suppose that $\hgt(I_{n-1}(\var)) \ge 3$.  Since 2 is a non
zerodivisor we may assume that $\var_{11}$, the $(1,1)$ entry of $\var$,
does not lie in $\im I_1(\var)$.  Having rank $n-1$, the matrix $\var$ fits
into an exact sequence
$$ 0 \longrightarrow R \overset \psi \to \longrightarrow R^n \overset \var
\to \longrightarrow R^{n*}. $$
As $\hgt(I_1(\psi)) \ge \hgt(I_{n-1}(\var)) \ge 3$, the complex
$$ F.: 0 \longrightarrow R \overset \psi \to \longrightarrow R^n \overset
\var \to \longrightarrow R^{n*} \overset \psi^* \to \longrightarrow R^*$$
is exact by the Buchsbaum-Eisenbud acyclicity criterion, see Buchsbaum and
Eisenbud [BE1, Theorem].  Thus $I_1(\psi) = I_1(\psi^*) \subset
I_1(\var)$.  Furthermore $I_1(\psi^*)$ is a Gorenstein ideal of height 3,
and hence by Buchsbaum and Eisenbud [BE2, Theorem 2.1], there is an exact sequence
$$ G.: 0 \longrightarrow R \overset \psi \to \longrightarrow R^n \overset
\chi \to \longrightarrow R^{n*} \overset \psi^* \to \longrightarrow R^*$$
with $\chi$ alternating.
 
The identity map on $\text{Ker}(\psi^*)$ lifts to a morphism of complexes
$\a.: F. \longrightarrow G.$ where $\a_0 = id$ and $\a_1 = id$.  Notice
that $\a_3$ is multiplication by some $u \in R$.  Thus, since $\a^*.: G.^*
\ra F.^*$ is a morphism of acyclic complexes of free modules, $\a.^*$ is
homotopic to multiplication by $u$.  Consequently $\a.$ has the same
property.  It follows that
$$ \var \equiv u \chi \text{ mod } (I_1(\chi)I_1(\var) + I_1(\chi)
I_1(\psi)),$$
hence
$$ \var \equiv u \chi \text{ mod } I_1(\var)^2.$$
But this is impossible because $\chi_{11} = 0$, whereas $\var_{11} \not\in
\im I_1(\var).$
\qed \enddemo

%%%%%%%%%%%%%%%%%%%%%%%%%%%%%%%%%%%%%%%%%%%%%%%%%%%%%%%%%%%%%%%%%%%%%%

\vskip1truein

\centerline{\bf References}
\bigskip
\refstyle{A}


\Refs\nofrills{}
\widestnumber\key{HH12}

\ref\key{AFH}
\by L. Avramov, H.-B. Foxby, and B. Herzog
\paper Structure of local homomorphisms
\jour J. Alg.
\vol 164
\yr 1994
\pages 124--145
\endref

\ref\key{B1}\by W. Bruns
\paper The Eisenbud-Evans Principal Ideal Theorem and determinantal
ideals
\jour Proc. Amer. Math. Soc.
\vol 83
\yr 1981
\pages 19--24
\endref








\ref
\key{BE1}
\manyby D. Buchsbaum and D. Eisenbud
\paper What makes a complex exact?
\jour J. Alg.
\vol 25
\yr 1973
\pages 259--268
\endref

\ref
\key{BE2}
\bysame
\paper Algebra structures for finite free resolutions and some
structure theorems for ideals of codimension 3
\jour Amer. J. Math.
\vol 99
\yr 1977
\pages 447--485
\endref



\ref
\key{CEU}
\by M. Chardin, D. Eisenbud and B. Ulrich
\paper Hilbert functions, residual intersections, and residually $S_2$ ideals
\paperinfo preprint 1998
\endref

\ref
\key{CN}
\by R.C. Cowsik and M.V. Nori
\paper On the fibres of blowing up
\jour J. Indian Math. Soc. 
\vol 40
\yr 1976
\pages 217--222
\endref

\ref
\key{EN}
\by J. Eagon and D. G. Northcott
\paper Ideals defined by matrices and a certain complex
associated to them.
\jour Proc. Royal Soc.
\vol 269
\yr 1962
\pages 188--204
\endref


\ref\key{EHU1}
\manyby D. Eisenbud, C. Huneke, and B. Ulrich
\paper A simple proof of some generalized principal ideal theorems
\jour Proc. Amer. Math. Soc.
\vol 129
\yr 2001
\pages 2535--254
\endref

\ref\key{EHU2}
\bysame
\paper What is the Rees algebra of a module?
\paperinfo to appear, Proc. Amer. Math. Soc.
\endref

\ref\key{EHU3}
\bysame
\paper Order ideals and a generalized Krull height theorem
\paperinfo preprint
\endref




\ref
\key {EG1}
\manyby Evans, E.G. and Griffith P.
\paper The syzygy problem
\jour Annals of Math.
\vol 114
\yr 1981
\pages 323--333
\endref


\ref
\key {EG2}
\bysame
\paper Syzygies
\jour London Math. Soc. Lecture Notes
\vol  106
\publ Cambridge Univ. Press, 1985
%\yr 1985
\endref



\ref
\key{F}
\by G. Faltings
\paper Ein Kriterium f\"ur vollst\"andige Durchschnitte
\jour Invent. Math.
\vol 62
\yr 1981
\pages 393--401
\endref

\ref
\key{FUV}
\by H. Flenner, B. Ulrich, and W. Vogel
\paper On limits of joins of maximal dimension
\jour Math. Ann.
\vol 308
\yr 1997
\pages 91--318
\endref

\ref
\key{G}
\by A. Grothendieck
\paper \'El\'ements de g\'eom\'etrie alg\'ebrique, III, \'Etude cohomologique
des faisceaux coh\'erents
\jour Institut des Hautes \'Etudes Scientifiques
\yr 1961
\vol 11
\endref

\ref\key{Ha1}
\manyby R. Hartshorne
\paper Varieties of small codimension in projective space
\jour Bull. Amer. Math. Soc.
\vol 80
\yr 1974
\pages 1017--1032
\endref

\ref\key{Ha2}
\bysame
\paper Complete intersections and connectedness
\jour Amer. J. Math.
\vol 84
\yr 1962
\pages 497--508
\endref





\ref
\key{Hu}
\by C. Huneke
\paper Criteria for complete intersections 
\jour  J. London Math. Soc. (2) 
\yr 1985
\vol 32
\pages 19--30
\endref








\ref\key {Ma}
\by H. Matsumura
\book Commutative Ring Theory
\publ Cambridge Univ. Press
\yr 1986
\endref

\ref\key{NR}
\by J. Eagon and D.G. Northcott
\paper Ideals defined by matrices and a certain complex associated with them
\jour Proc. Roy. Soc. Ser. A
\vol 269
\yr 1962
\pages 188--204
\endref

\ref
\key {R}
\by D. Rees
\paper Reduction of modules
\jour Math. Proc. Camb. Phil. Soc.
\vol 101
\yr 1987
\pages 431--449
\endref




\ref\key{S}
\by J.-P. Serre
\book Alg\`ebre locale, multiplicit\'es
\bookinfo Springer Lect. Notes in Math.
\vol  11
\yr  1958
\endref


\ref\key{SUV}
\by A. Simis, B. Ulrich and W. Vasconcelos
\paper Rees algebras of modules
\paperinfo preprint
\endref


\endRefs
\enddocument